\documentclass[letter, 10pt, conference]{ieeeconf}  

\IEEEoverridecommandlockouts                              
\usepackage{amsmath,graphicx}
\usepackage{amssymb}
\usepackage[dvips]{epsfig}
\usepackage{subfigure}
\usepackage{multirow}
\DeclareMathOperator*{\argmin}{arg\,min}

\title{\LARGE \bf 
Lyapunov-based Stochastic Nonlinear Model Predictive Control: \\ Shaping the State Probability Density Functions }

\author{Edward A. Buehler$^1$, Joel A. Paulson$^{1,2}$, Ali Akhavan$^1$, and Ali Mesbah$^{1,\dagger}$ 
\thanks{$^1$Department of Chemical and Biomolecular Engineering, University of California, Berkeley, CA 94720, USA.}
\thanks{$^2$Department of Chemical Engineering, Massachusetts Institute of Technology, MA 02139, USA.}
\thanks{$^{\dagger}$Corresponding author: {\tt\small mesbah@berkeley.edu}.}
}

\newcounter{myAlgoCounter}
\usepackage{mathtools}
 \usepackage{ssCorrections}
 \xdefinecolor{colSS}{rgb}{1,0,0} 

\def \r {\mathbb{R}}
\def \x {\mathbf{x}}
\def \u {\mathbf{u}}
\def \w {\mathbf{w}}
\def \f {\mathbf{f}}
\def \g {\mathbf{g}}
\def \h {\mathbf{h}}
\def \k {\mathbf{k}}
\def \Q {\mathbf{Q}}
\def \d {d}  
\def \V {\mathbf{V}}
\def \pc {\mathbf{p}}
\def \tr {\mathsf{Tr}}
\def \L {\mathcal{L}}
\def \U {\mathbb{U}}

\def \p {\mathsf{P}}
\def \pr {\mathbf{Pr}}

\begin{document}

\maketitle
\thispagestyle{empty}
\pagestyle{empty}

\begin{abstract}
Stochastic uncertainties in complex dynamical systems lead to variability of system states, which can in turn degrade the closed-loop performance. This paper presents a stochastic model predictive control approach for a class of nonlinear systems with unbounded stochastic uncertainties. The control approach aims to shape probability density function of the stochastic states, while satisfying input and joint state chance constraints. Closed-loop stability is ensured by designing a stability constraint in terms of a stochastic control Lyapunov function, which explicitly characterizes stability in a probabilistic sense. The Fokker-Planck equation is used for describing the dynamic evolution of the states' probability density functions. Complete characterization of probability density functions using the Fokker-Planck equation allows for shaping the states' density functions as well as direct computation of joint state chance constraints. The closed-loop performance of the stochastic control approach is demonstrated using a continuous stirred-tank reactor.                
\end{abstract}

\section{Introduction}
\label{sec:Intro}

The need to account for system uncertainties in model predictive control (MPC) of complex dynamical systems has led to extensive investigation of robust MPC approaches (e.g., see \cite{may14} and the references therein). The majority of work on robust MPC considers bounded, deterministic uncertainty descriptions with the goal to design MPC control laws that are robust to \textit{worst-case} system uncertainties. The deterministic robust MPC approaches may, however, result in conservative closed-loop control performance, as worst-case system uncertainties are likely to have a small probability of occurrence \cite{bem99}. This consideration has recently motivated the development of stochastic MPC (SMPC) approaches that directly use probabilistic descriptions of the stochastic system uncertainties (i.e., parametric uncertainties, uncertain initial conditions, and exogenous disturbances). In particular, SMPC approaches allow for defining chance constraints in the stochastic optimal control problem to systematically seek tradeoffs between the control performance and robustness to system uncertainties.

The formulation of a SMPC approach largely depends on the complexity of system dynamics, properties of stochastic uncertainties, and solution method for the stochastic programming problem. SMPC approaches have been proposed for linear systems with multiplicative noise \cite{can09,pri09} and additive noise \cite{old10,hok12,far13,pau15}. The latter approaches mainly use affine parameterizations of control inputs for finite-horizon linear quadratic problems to transform the stochastic programming problem into a deterministic one. Randomized algorithms have also been used to develop SMPC approaches for linear systems \cite{bat04,cal13}. Recently, a SMPC approach has been proposed for nonlinear systems with time-invariant probabilistic uncertainties using the generalized polynomial chaos framework \cite{mes14} (also see \cite{hes06} for SMPC for nonlinear systems in the absence of input constraints). Generally, the characteristics of stochastic uncertainties (e.g., boundedness, additive/multiplicative, and time-varying/time-invariant nature of stochastic uncertainties) have important implications for closed-loop stability and recursive feasibility properties of SMPC approaches with input constraints and chance constraints. In addition, the complexity of system dynamics largely affects the computational complexity of probabilistic uncertainty propagation (through system dynamics) as well as chance constraint handling.

This paper presents a stochastic nonlinear MPC (SNMPC) approach for a class of nonlinear systems with probabilistic uncertain initial conditions and unbounded stochastic disturbances. The proposed SNMPC approach includes input constraints and joint state chance constraints. To ensure closed-loop stability of the SNMPC approach, a stochastic Lyapunov-based feedback control law that explicitly characterizes stability in a probabilistic sense is used (e.g., \cite{den00,mah12}). The Lyapunov-based feedback control law allows for designing a \textit{stability constraint} in terms of a stochastic control Lyapunov function, which guarantees that the origin of the closed-loop system is asymptotically stable in probability (Section~\ref{sec:Prel}).

The Lyapunov-based SNMPC approach is intended to shape probability density functions (PDFs) of the stochastic state variables. This necessitates characterizing the complete PDFs of states. The Fokker-Planck equation \cite{ris84} is used to describe the dynamic evolution of the (multivariate) PDFs associated with the stochastic nonlinear system (Section~\ref{sec:FP}). Complete characterization of the states' PDFs also allows for direct computation of joint state chance constraints without approximation. This work uses the Hellinger distance \cite{gib02} to quantify the similarity between the predicted (multivariate) PDFs of states and user-specified reference PDFs for shaping the probability density functions (Section~\ref{sec:hellinger}). Note that control of PDFs of system states/outputs \cite{jum92,wan99,for04} and MPC of stochastic systems using the Fokker-Planck equation \cite{bor13,fle14} have been reported in the literature. What distinguishes this work is the generic formulation of the Lyapunov-based SNMPC approach in terms of PDF shaping as well as input and joint state chance constraints handling, while ensuring closed-loop stability (Section~\ref{sec:SNMPC_formulation}). The presented Lyapunov-based SNMPC approach is demonstrated for stochastic optimal control of a continuous stirred-tank reactor in the presence of stochastic uncertainties (Section~\ref{sec:CSTR}).

\section{Preliminaries}
\label{sec:Prel}

\noindent \textit{\textbf{Notation} }

\vspace{2mm}

Throughout this paper, boldface symbols (e.g., $\x$) denote vectors and subscripts denote vector elements (e.g., $x_i$). $\r^n$ denotes the $n$-dimensional Euclidean space with $\r_{+}=[0, \; \infty)$. The transpose of a vector or a matrix will be denoted by superscript $\top$. $\tr\{\cdot\}$ denotes the trace operator on a square matrix. For a vector $\x \in \r^n$, $\left\| \x \right\|$ denotes the Euclidean norm of $\x$, and $\left\| \x \right\|^2_{\Q}$ denotes the weighted norm of $\x$ defined by $\left\| \x \right\|^2_{\Q} = \x^{\top}\Q\x$ with $\Q$ being a positive definite symmetric matrix. $\p_{\x}$ denotes the (multivariate) probability density function of $\x \in \r^n$. $(\Omega,\mathcal{F},\mathcal{P})$ denotes a probability space defined by the sample space $\Omega$, $\sigma$-algebra $\mathcal{F}$, and probability measure $\mathcal{P}$ on $\Omega$. $\pr\{\cdot\}$ denotes the probability of satisfaction of an expression. $\L_\f\mathcal{X}$ denotes the Lie derivative of a scalar function $\mathcal{X}(\cdot)$ with respect to a vector function $\f(\cdot)$. A continuous function $\V: \r^n \rightarrow \r$ is said to be $C^k$ if it is $k$-times differentiable. A continuous function $\alpha: \r_{+} \rightarrow \r_{+}$ is said to belong to class $\mathcal{K}$ when it is strictly increasing and $\alpha(0)=0$. The function $\alpha$ is said to belong to class $\mathcal{K}_{\infty}$ when $\alpha \in \mathcal{K}$ and $\alpha(a) \rightarrow \infty$ as $a \rightarrow \infty$. 

\vspace{2mm}

\noindent \textit{\textbf{System description} }        

\vspace{2mm}

Consider a class of stochastic nonlinear systems described by the stochastic differential equation (SDE) 
\begin{equation} \label{e_sys}
\begin{array}{ll}
 \d \x(t)  & = \f(\x(t))\d t + \g(\x(t))\u(t)\d t + \h(\x(t))\d \w(t) \\
 \x(t_0) & \sim \p_{\x_0},
\end{array}
\end{equation}  
where $\x(t) \in \r^{n}$ denotes the stochastic state variables with the known initial multivariate PDF $\p_{\x_0}$; $\u(t) \in \r^{m}$ denotes the system inputs; $\w(t)$ denotes a $q$-dimensional standard Wiener process (i.e., stochastic disturbances) defined on the probability space $(\Omega,\mathcal{F},\mathcal{P})$; and $\f:\r^n \rightarrow \r^n$, $\g:\r^n \rightarrow \r^{n\times m}$, and $\h:\r^n \rightarrow \r^{n\times q}$ denote the Borel measurable functions that describe the system dynamics. The functions $\f$, $\g$, and $\h$ are assumed to be locally bounded and locally Lipschitz continuous in $\x(t), \; \forall t \in \r_{+}$, and $\f(0)=0$ (i.e., the origin is the steady-state point of the unforced and undisturbed system). The latter conditions ensure uniqueness and local existence of solutions to the SDE~\eqref{e_sys} \cite{kar91}. The system inputs $\u(t)$ in~\eqref{e_sys} are constrained to lie in a nonempty convex set $\U \subseteq \r^{m}$ defined by
\begin{equation} \label{e_ic}
\U \coloneqq \{ \u(t) \in \r^{m} \; | \; \u_{min} \leq \u(t) \leq \u_{max} \},
\end{equation}        
where $\u_{min} \in \r^{m}$ and $\u_{max} \in \r^{m}$ denote the lower and upper bounds on $\u$, respectively. In addition, the stochastic system states $\x(t)$ should satisfy hard inequality constraints  
\begin{equation} \label{e_sc}
\k(\x(t)) \leq 0,
\end{equation}
where $\k:\r^n \rightarrow \r^p$ denotes (possibly) nonlinear functions that describe the state constraints, and $\k(0)=0$.    


Note that the stochasticity of system~\eqref{e_sys} arises from the probabilistic nature of uncertain initial states $\x(t_0)$ (described by $\p_{\x_0}$) and the stochastic disturbances $\w$. In \eqref{e_sys}, the terms $\f(\x(t)) + \g(\x(t))\u(t)$ and $\h(\x(t))$ correspond to the drift and diffusion terms in the \textit{Ito stochastic process}, respectively \cite{gar09}. Any general (deterministic) nonlinear model can be represented in terms of the control-affine deterministic drift term $\f(\x(t)) + \g(\x(t))\u(t)$ \cite{nij90}.

\vspace{2mm}

\noindent \textit{\textbf{Stochastic optimal control with state chance constraints} }        

\vspace{2mm}

This paper investigates stochastic MPC of the nonlinear system~\eqref{e_sys} such that stability of the closed-loop system is guaranteed. The proposed SNMPC approach should allow for shaping the multivariate PDF $\p_\x(t)$ of system states in an optimal manner, while the system inputs $\u(t)$ lie in the set $\U$. In addition, the stochastic optimal control approach should ensure satisfaction of the (possibly nonlinear) state constraints~\eqref{e_sc} with at least probability $\beta$ in the presence of system stochasticity. This requires incorporating \textit{joint} chance constraints of the form 
\begin{equation} \label{e_cc}
\pr\{\k(\x(t)) \leq 0\} \geq \beta
\end{equation}
into the stochastic optimal control problem. This paper considers \textit{receding-horizon} implementation of the SNMPC approach in a full state feedback control scheme, where the PDF $\p_{\x}(t_k)$ is assumed to be known at every measurement sampling time instant $t_k$.          

The key challenges that will be addressed in this paper for solving the above described stochastic optimal control problem are: (i) describing the dynamic evolution of the multivariate PDF $\p_{\x}(t)$ associated with the SDE~\eqref{e_sys}, (ii) converting the joint chance constraint~\eqref{e_cc} to computationally tractable expressions, and (iii) designing the SNMPC control law such that it ensures closed-loop stability of the stochastic system. Next, the main result of stochastic Lyapunov stability (see \cite{den00,den01}) that will be used for designing a Lyapunov-based SNMPC approach is summarized.      

\vspace{2mm}

\noindent \textit{\textbf{Lyapunov-based controllers} }        

\vspace{2mm}

For stochastic nonlinear systems, Lyapunov-based stabilizing control laws that explicitly characterize region of attraction of the closed-loop system in a probabilistic sense have been proposed (e.g., see \cite{den00,cha99,ste03,mah12}, and the references therein). In this paper, a Lyapunov-based control law is designed for the stochastic optimal control problem presented above. It is assumed that there exists a nonlinear feedback control law $\u(t)=\pc(\x), \; \forall \x \in \mathcal{X} \subseteq \r^n$, where $\pc:\r^n \rightarrow \r^m$ denotes a nonlinear function and $\mathcal{X}$ denotes a compact set that contains the origin $\x=0$. The feedback control law $\pc(\x)$ is intended to make the closed-loop system asymptotically stable (in probability) about the origin, while the input and state constraints~(i.e., \eqref{e_ic} and \eqref{e_cc}) are satisfied. According to the converse Lyapunov theorem \cite{kha96}, the existence of the feedback control law $\pc(\x)$ implies the existence of a \textit{stochastic control Lyapunov function} $\V(\x)$ that is defined as in Thm.~1.

 \textbf{Theorem 1 (Asymptotic stability in probability \cite{den00}):} Consider the stochastic nonlinear system~\eqref{e_sys} and assume that there exists a $C^2$-function $\V:\r^n \rightarrow \r_{+}$, class $\mathcal{K}_{\infty}$ functions $\alpha_1$ and $\alpha_2$, and a class $\mathcal{K}$ function $\alpha_3$, such that $\forall \x\in \mathcal{X},\; \forall t\geq 0$  
\begin{equation*}
\alpha_1(\left|\x\right|) \leq \V(\x) \leq \alpha_2(\left|\x\right|) ,
\end{equation*}        
\begin{equation*}
\begin{split}
\L_\f \V(\x) & + \L_\g \V(\x)\u(t)|_{\u(t)=\pc(\x)} +  \\ 
& \frac{1}{2}\tr\{\h(\x)^\top\frac{\partial^2 \V}{\partial \x^2}\h(\x)\} \leq \alpha_3(\left|\x\right|).
\end{split}
\end{equation*}  
Then the stochastic control Lyapunov function $\V(\x)$ ensures that the origin is asymptotically stable in probability. $ \; \; \quad  \blacksquare$   

Thm.~1 indicates that the stochastic Lyapunov-based control techniques allow for defining feedback control laws that will lead to     
\begin{equation} \label{e_lyp}
\begin{split}
 \L_\f \V(\x) & + \L_\g \V(\x)\u(t)|_{\u(t)=\pc(\x)}  + \frac{1}{2}\tr\{\h(\x)^\top\frac{\partial^2 \V}{\partial \x^2}\h(\x)\} \\ 
&  + \gamma \V(\x) \leq 0,  \quad \quad \quad \quad \forall \x\in \Pi,
\end{split}
\end{equation}
where the set $\Pi$ is defined by
\begin{equation*}
\Pi \coloneqq \sup_{c\in \r} \; \{ \x \in \r^{n} \; | \; \u \in\U, \; \x \in \mathcal{X}, \; \V(\x)\leq c \};
\end{equation*} 
and $\gamma > 0$ is a constant. Next, the stochastic control Lyapunov function $\V(\x)$ will be used for designing a control law for the stochastic optimal control problem such that closed-loop stability of the proposed SNMPC approach is guaranteed.    

\section{Stochastic Nonlinear Model \\ Predictive Control}
\label{sec:SNMPC}

This section presents the formulation of the Lyapunov-based SNMPC approach with joint state chance constraints. The propagation of probabilistic system uncertainties (i.e., uncertain initial states and stochastic disturbances) through system dynamics is described by the \textit{Fokker-Planck equation}. The \textit{Hellinger distance} is used as a measure of similarity of multivariate PDFs to formulate the objective function of the stochastic optimal control problem for PDF shaping.    

\subsection{Fokker-Planck Equation for Uncertainty Propagation}
\label{sec:FP}

The Fokker-Planck (FP) equation describes the dynamic evolution of the PDF of stochastic states $\x(t)$ in the uncertain nonlinear system~\eqref{e_sys} \cite{ris84}. The FP equation readily characterizes the complete multivariate PDF $\p_{\x}(t)$ arisen from the stochastic system uncertainties in initial states $\x(t_0)$ and disturbances $\w$. This is in contrast to uncertainty propagation techniques that describe merely certain \textit{statistics} of the PDFs (e.g., see \cite{mes14} and the references therein). Characterizing the complete PDF of states using the FP equation enables the proposed SNMPC approach to: (i) shape the PDF $\p_{\x}$ with respect to any desired (multivariate) PDF, and (ii) compute chance constraints of any complexity directly without conservative approximations.

The FP equation associated with the SDE~\eqref{e_sys} that describes the evolution of the multivariate PDF $\p_{\x}(t)$ is defined by
\begin{equation} \label{e_FP}
\begin{split}
\frac{\partial \p_{\x}}{\partial t} &+ \sum_{i=1}^{n} \frac{\partial }{\partial x_i} \bigg ( \big( f_i(\x) + g_i(\x)\u \big ) \p_{\x} \bigg)\\
& - \frac{1}{2} \sum_{i=1}^{n} \sum_{j=1}^{n} \frac{\partial^2 }{\partial x_i \partial x_j} \bigg (D_{ij}(\x)  \p_{\x} \bigg) = 0
\end{split}
\end{equation} 
with the initial condition
\begin{equation*} 
\p_{\x}(0) = \p_{\x_0}, 
\end{equation*} 
where $\mathbf{D}=\h(\x)\h(\x)^\top$ (i.e., the diffusion matrix). The FP equation~\eqref{e_FP} is a parabolic partial differential equation, whose solution should be nonnegative and satisfy    
\begin{equation*} 
\int_{\Omega_{\x}} \p_{\x}(t) \d \x=1, \quad \forall t \geq 0 . 
\end{equation*}
The existence and uniqueness of a solution to~\eqref{e_FP} under mild assumptions have been established (see \cite{ris84,jor98}). Note that the FP equation~\eqref{e_FP} can be used to compute the univariate PDFs $\p_{x_i}$ for every state $x_i$ and joint PDFs for any combination of stochastic states.

Solving the FP equation is generally challenging for nonlinear systems, in particular systems with high state dimension \cite{ris84}. Various numerical methods such as finite difference and finite element methods have been used to solve the FP equation for nonlinear systems (e.g., see \cite{kum06} and the references therein). In this work, finite volume method with first order upwind interpolation scheme is used to solve~\eqref{e_FP} \cite{lev02}. The finite volume method allows for effectively dealing with the convective nature of the FP equation (due to the drift term), as well as suppressing numerical instability problems.

\subsection{Metric for Similarity of Probability Density Functions}
\label{sec:hellinger}

The SNMPC approach aims to shape the PDF $\p_{\x}$ according to a predetermined (arbitrarily-shaped) PDF. Hence, a metric is required to quantify the similarity between the predicted and the reference PDFs at each time instant $t$. 

This paper uses the Bhattacharyya coefficient \cite{kai67}, a measure closely related to the \textit{Bayes error} \cite{and03}, to establish a measure for the similarity of PDFs. The Bhattacharyya coefficient quantifies the degree of overlap between two (multivariate) PDFs, and is defined by
\begin{equation} \label{e_bha}
\mathfrak{B}(\x)  \coloneqq \int_{\Omega_{\x}} \sqrt{\p_{\x}(t) \p_{\x}^{\textbf{ref}} } \d \x,  
\end{equation}
where $\p_{\x}^{\textbf{ref}}$ denotes a reference PDF. The Bhattacharyya coefficient will be larger when the overlap between the PDFs is larger. $\mathfrak{B}(\x)=0$ if the PDFs do not overlap, whereas $\mathfrak{B}(\x)=1$ when the PDFs are identical. Explicit forms of the Bhattacharyya coefficient for various PDFs are given in \cite{kai67,djo90}.    

The Bhattacharyya coefficient~\eqref{e_bha} is used to define a metric for the similarity between $\p_{\x}(t)$ and $\p_{\x}^{\textbf{ref}}$ in the objective function of the stochastic optimal control problem. The metric, which is known as the Hellinger distance \cite{gib02}, is defined by
\begin{equation} \label{e_hel}
\Delta(\x) \coloneqq \sqrt{1-\mathfrak{B}(\x)}. 
\end{equation}
Note that the metric~\eqref{e_hel} is near optimal due to its relation to the Bayes error, and can be used for arbitrary PDFs \cite{com00}.

\subsection{Formulation of the Lyapunov-based SNMPC Approach}
\label{sec:SNMPC_formulation}

The Fokker-Planck equation~\eqref{e_FP} and the Hellinger distance~\eqref{e_hel} are now used to formulate the Lyapunov-based SNMPC problem for the stochastic nonlinear system~\eqref{e_sys}.   

 \textbf{Problem 1}  \textbf{(Lyapunov-based SNMPC with input and joint state chance constraints):} Suppose that the PDF $\p_{\x}(t_k)$ is known at every sampling time instant $t_k$.\footnote{Note that this work considers a full state feedback control scheme. The PDFs $\p_{\x}(t_k)$ arise from measurement errors.} The stochastic optimal control problem at each time instant $t_k$ is stated as  
\begin{align*}
\u^\ast(\p_{\x}(t_k)) \coloneqq \underset{\u }{\argmin} \int_{0}^{T_p} \bigg (\Delta (\bar{\x}(t')) + \left\| \u(t')  \right\|_{\mathbf{R}}^2  \bigg) d t'
\end{align*} 
\begin{align} \label{e_P1}
\begin{array}{lll}
\text{s.t.:} & {\displaystyle \frac{\partial \p_{\bar{\x}}}{\partial t} }+ \smashoperator[r]{\sum_{i=1}^{n}} \frac{\partial }{\partial \bar{x}_i} \bigg ( \big( f_i(\bar{\x}) + g_i(\bar{\x})\u \big ) \p_{\bar{\x}} \bigg) \\
& - \frac{1}{2} \smashoperator[r]{\sum_{i=1}^{n}} \sum_{j=1}^{n} \frac{\partial^2 }{\partial \bar{x}_i \partial \bar{x}_j} \bigg (D_{ij}(\bar{\x})  \p_{\bar{\x}} \bigg) = 0, & \forall t \in [0,\; T_p]  \\ \\ 
&   \L_\f \V(\bar{\x}) + \L_\g \V(\bar{\x})\u(t)  +    \\
&  \frac{1}{2}\tr\{\h(\bar{\x})^\top\frac{\partial^2 \V}{\partial \bar{\x}^2}\h(\bar{\x})\} + \gamma \V(\bar{\x}) \leq 0, & \forall t \in [0,\; T_p]   \\  \\
& \pr\{\k(\bar{\x}(t)) \leq 0\} \geq \beta, & \forall t \in [0,\; T_p] \\ \\
& \u(t) \in \U, & \forall t \in [0,\; T_c] \\ \\
& \p_{\bar{\x}}(0) = \p_{\x}(t_k),
\end{array}
\end{align}
where $\u^\ast$ denotes the optimal inputs (i.e., control policy) over the control horizon $[0,\; T_c]$; $\bar{\x}$ denotes the states predicted by the system model; $T_p$ denotes the prediction horizon; and $\mathbf{R}$ denotes a strictly positive definite matrix. The closed-loop stability of the SNMPC approach is ensured in a probabilistic sense by incorporating the stability constraint (defined in terms of the stochastic control Lyapunov function $\V(\bar{\x})$) into \eqref{e_P1}. The asymptotic stability properties of the closed-loop system are characterized in Thm.~1 and \eqref{e_lyp}. The proposed stochastic control approach possesses the stability properties of Lyapunov-based controllers when applied in a sample-and-hold fashion (e.g., see \cite{mah12}). Note that the stochastic optimal control problem~\eqref{e_P1} is implemented in a receding-horizon mode, where merely the optimal inputs $\u^{\ast}(0)$ are applied to the stochastic nonlinear system~\eqref{e_sys} at every sampling time instant $t_k$.

The Lyapunov-based SNMPC approach in Problem 1 allows for shaping the multivariate PDF of states, while satisfying input constraints and joint chance constraints imposed on the stochastic states. This is due to using the FP equation for probabilistic uncertainty propagation, as the FP equation enables explicit characterization of the states' PDFs. The objective function in~\eqref{e_P1} is stated in its most general form in terms of the Hellinger distance to quantify the difference between the multivariate PDFs $\p_{\x}(t)$ and $\p_{\x}^{\textbf{ref}}$ over the prediction horizon $t \in [0,\; T_p]$.\footnote{When the control objective is to achieve desired PDFs for individual states, the objective function in~\eqref{e_P1} can be defined in terms of weighted sum of Hellinger distances pertaining to the univariate PDF of states.} The objective function can be simplified by considering only certain statistics (e.g., the expected value and variance) of the PDFs. Furthermore, the joint state chance constraint in~\eqref{e_P1} can be readily computed through explicit integration without any approximation since the knowledge of the full multivariate PDF is available. When merely individual state chance constraints are considered in~\eqref{e_P1}, the FP equation can be adapted to compute only the univariate PDFs for the respective states. This will reduce the computational complexity of the stochastic optimal control problem. Next, the application of the Lyapunov-based SNMPC approach to a continuous stirred-tank reactor is investigated.            

\section{Case Study: Stochastic Optimal Control of a Continuous Stirred-Tank Reactor}
\label{sec:CSTR}

Consider a continuous stirred-tank reactor (CSTR), in which the exothermic reaction $A \xrightarrow{k} B$ occurs. The system dynamics are described by
\begin{align*} 
 \d C_A  = & \bigg (  \frac{F}{V}(C_{A0}-C_A)-k_0e^{\frac{-E}{RT}}C_A \bigg ) \d t \\
& +  \sigma_{C_A} \d w_{C_A}(t), \quad C_A(0) \sim \mathcal{B}(0,2,320,320)  \\ \\
 \d T  = & \bigg (  \frac{F}{V}(T_{0}-T) + \frac{\delta H}{\rho c_p}k_0e^{\frac{-E}{RT}}C_A + \frac{Q}{\rho c_p V} \bigg )\d t, \\
&  T(0) = 315.0,  
\end{align*}  
where $C_A$ denotes the concentration of species $A$ (kmol/m$^3$); $T$ denotes the reactor temperature (K); $C_{A0}$ (with the mean value of $0.702$ kmol/m$^3$) and $T_0$ denote the concentration of species $A$ and the temperature in the inlet reactor stream, respectively; $F$ denotes the inlet flow rate (m$^3$/min); $V$ denotes the reactor volume (m$^3$); $\mathcal{B}$ denotes the four-parameter beta distribution; $Q$ denotes the heat removed from the reactor (kJ/min); and $w_{C_A}(t)$ denotes a standard Wiener process acting on $C_A$. The model parameters are listed in Table~\ref{T1}. Note that the CSTR under study is a stochastic system due to the probabilistic uncertainties in the initial concentration $C_A(0)$, as well as the stochastic disturbances $w_{C_A}(t)$. The system inputs that can be manipulated for control are $C_{A0}$ and $Q$. It is assumed that the temperature $T$ and the PDF of $C_A$ (which is of beta-distribution type) are measured at every sampling time instant $t_k$ (i.e., $2$ min). The process is run for $30$ min.

The Lyapunov-based SNMPC approach presented in Problem 1 is applied to stabilize the CSTR around the steady-state point ($\bar{C}_A = 0.57$ kmol/m$^3$, $\bar{T}=317$ K), while shaping the probability density of $C_A$ to take the form of the Normal Distribution $\mathcal{N}(\bar{C}_A,4\times10^{-4})$ (i.e., $\p_{C_A}^{\textbf{ref}}=\mathcal{N}(\bar{C}_A,4\times 10^{-4})$ in \eqref{e_bha}). The objective function in~\eqref{e_P1} is defined as
\begin{align} \label{e_obj} 
\int_{0}^{T_p} \bigg (\Delta (C_A(t')) + \left\| \mathbb{E}[T(t')] - \bar{T}  \right\|^2  \bigg) d t',
\end{align}
where $\mathbb{E}[\cdot]$ denotes the expected value. The prediction horizon and the control horizon are selected to be $T_p=30$ min and $T_c=20$ min, respectively. Note that the objective function~\eqref{e_P1} is defined such that the SNMPC approach minimizes the difference between the PDF of $C_A$ and the reference PDF $\p_{C_A}^{\textbf{ref}}$ over the prediction horizon $[0,\;T_p]$. To design the stability constraint in~\eqref{e_P1}, the quadratic Lyapunov function $\V(\x)$ is defined as $\V(\bar{\x})=\bar{\x}^\top\mathbf{P}\bar{\x}$, where $\bar{\x}=[C_A-\bar{C}_A \; \; T-\bar{T}]^\top$ and $\mathbf{P}=\left[\begin{array}{cc} 3.18 & 0.93 \\ 0.93 & 0.58 \end{array} \right]$. The matrix $\mathbf{P}$ is obtained by solving the Lyapunov equation for the (nominal) linearized system dynamics around the steady-state point.   

\begin{table}[t!]
\caption{CSTR model parameters.}
\begin{center}
\begin{tabular}{llll}
 \hline  \\
$V$ & $0.1$~m$^3$ &  $\delta H$ & $4.78\times10^5$~kJ/kmol \\
$F$ & $100\times10^{-3}$~m$^{3}$/min & $k_0$ & $72\times10^9$~min$^{-1}$ \\
$T_0$ & $315.0$~K & $c_p$ & $0.239$~kJ/kgK  \\
$E$ & $8.314\times10^4$~kJ/kmol & $\rho$ & $1000$~kg/m$^{3}$  \\
$R$ & $8.314$~kJ/kmol K & $\sigma_{C_A}$ & $0.32$  \\ [1.0ex]  \hline
\end{tabular}  
\end{center}
\label{T1}
\end{table}

\begin{figure}[b!] 
\centering
\includegraphics[width=250pt]{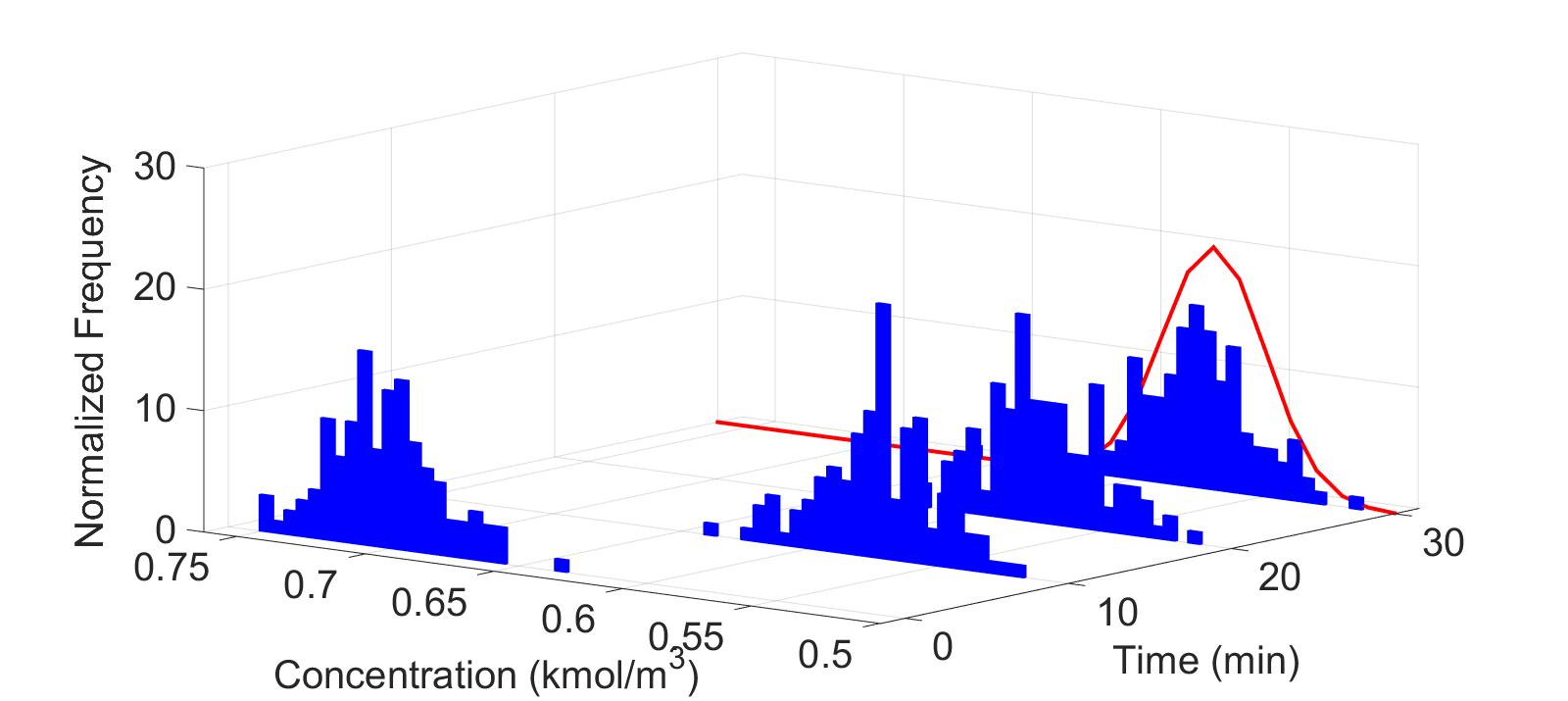}
\caption{Histograms of the concentration $C_A$ at various times based on $130$ closed-loop simulations with different disturbance realizations $w_{C_A}$. The SNMPC approach is intended to shape the probability distribution of $C_A$ in the course of the process according to the reference probability distribution (red solid distribution).}
\label{fig1}
\end{figure} 

In the stochastic optimal control problem~\eqref{e_P1}, the inputs to the system are constrained to lie in the ranges $0 \leq C_{A0} \leq 2$ kmol/m$^3$ and $|Q|\leq 10$ kJ/min. In addition, the concentration $C_A$ should remain above a threshold in the presence of probabilistic system uncertainties. Hence, the individual chance constraint 
\begin{align}  \label{e_cccs} 
\pr\{C_A(t) \leq 0.53\} \leq 0.05
\end{align}
is incorporated into~\eqref{e_P1}, suggesting that the state constraint $C_A(t) \leq 0.53$ should be satisfied with at least probability $95\%$ in a stochastic setting. Since the individual chance constraint and the univariate PDF shaping in the objective function have been defined merely in terms of the concentration $C_A$, the FP equation~\eqref{e_FP} is considered only for $C_A$ with the diffusion coefficient $D=0.001$. The FP equation is solved using the finite volume method with $200$ discretization points over the concentration support $[0,\;2]$. The above constrained nonlinear optimal control problem is solved using the MATLAB subroutine \texttt{fmincon}, where the set of model equations is integrated using the solver \texttt{ODE45}. The control inputs are discretized in a piecewise constant fashion in $5$ intervals over the control horizon.

\begin{figure}[t!] 
\centering
\includegraphics[width=250pt]{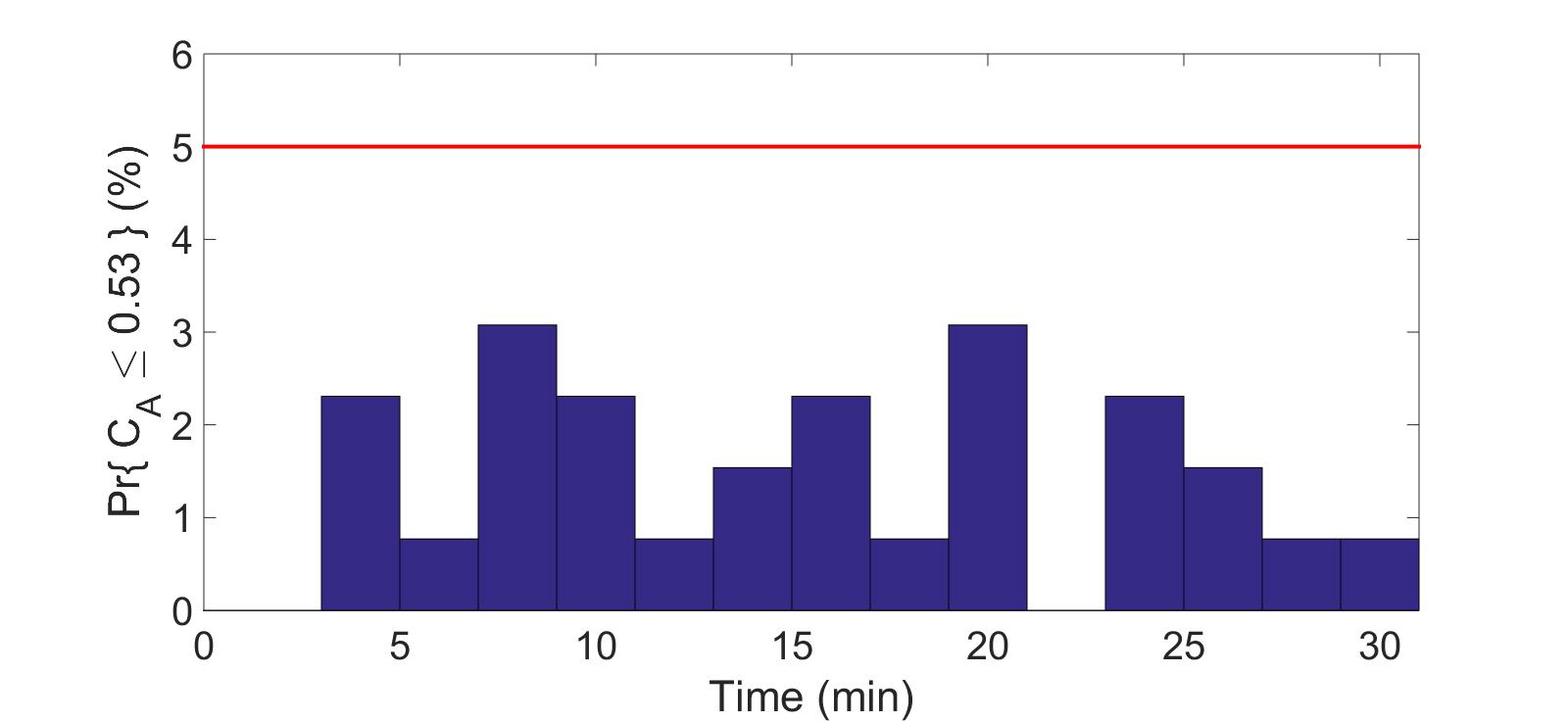}
\caption{Probability of state constraint violation at various times during the process. The probability of state constraint violation always remains below $5\%$ (red solid line) due to the state chance constraint~\eqref{e_cccs}.}
\label{fig2}
\end{figure}

To evaluate the performance of the SNMPC approach, Monte Carlo simulations of the closed-loop system are performed based on $130$ realizations of the Wiener process $w_{C_A}$  (the disturbance realizations are used for simulating the CSTR model to which
the optimal control inputs are applied at every sampling time instant $t_k$). Fig.~\ref{fig1} shows the evolution of the probability distribution of the concentration $C_A$ (i.e., true plant outputs) in the course of the process. The SNMPC approach stabilizes the stochastic CSTR system around the steady-state point $\bar{C}_A = 0.57$ kmol/m$^3$; the CSTR temperature is also stabilized around its steady-state value $\bar{T}=317$ K (not shown here). In addition, Fig.~\ref{fig1} shows that the SNMPC approach can effectively shape the PDFs of $C_A$ according to the reference PDF $\p_{C_A}^{\textbf{ref}}$. For instance, the mean and variance of the PDF of $C_A$ at time $30$ min are $0.576$ and $4.3\times10^{-4}$, respectively, which are aligned with those of the reference PDF $\p_{C_A}^{\textbf{ref}}=\mathcal{N}(0.57,4\times 10^{-4})$. The probability of violation of the state constraint over the process time is shown in Fig.~\ref{fig2}. The probability of constraint violation remains below the predetermined probability level $5\%$ (i.e., $\beta = 0.05$ in \eqref{e_cccs}) at all times during the process. Effective constraint handling is due to the state chance constraint~\eqref{e_cccs}, which ensures constraint satisfaction with a desired probability level in the presence of stochastic uncertainties. Note that simulation results (not shown here) revealed that by decreasing the desired probability of chance constraint satisfaction, the closed-loop control performance in terms of the PDF shaping can be further improved. This indicates the capability of the SNMPC approach to systematically seek tradeoffs between the closed-loop performance and robustness to system stochasticities.           

\section{Conclusions}
\label{sec:Conclusions}

A stochastic model predictive control approach is presented for a class of nonlinear systems with stochastic uncertainties. The closed-loop stability of the control approach is ensured by explicitly characterizing stability in a probabilistic sense using a stochastic control Lyapunov function. A key challenge in SMPC is efficient propagation of uncertainties through the system dynamics to fully characterize the probability distribution of the stochastic states. This paper uses the Fokker-Planck equation for uncertainty propagation, which allows for describing evolution of the probability distributions of the stochastic states for general probabilistic uncertainty descriptions. This work demonstrates that characterization of the complete probability distribution enables shaping the states' probability distributions with respect to arbitrarily-shaped reference distributions, as well as direct computation of chance constraints without any approximation.

%
%

\bibliographystyle{ieeetr}
\bibliography{Literature_list}

\end{document}